\documentclass[12pt,russian]{article}

\usepackage{amsfonts,amsmath,amsthm,amssymb,pb-diagram}

    \usepackage[T2A]{fontenc}
    \usepackage[cp1251]{inputenc}
    \usepackage[russian]{babel}

\numberwithin{equation}{section}

\theoremstyle{plain}
\newtheorem{theorem}{Теорема}[section]
\newtheorem{lemma}{Лемма}[section]

\newtheorem{statement}{Утверждение}[section]
\newtheorem{corollary}{Следствие}[section]
\theoremstyle{definition}
\newtheorem{definition}{Определение}[section]
\newtheorem{remark}{Замечание}[section]
\newtheorem{example}{Пример}[section]

\DeclareMathOperator{\id}{id}
\DeclareMathOperator{\card}{card}

\DeclareMathOperator*{\Lim}{Lim}

\newcommand{\abs}[1]{\left\vert#1\right\vert}
\newcommand{\set}[1]{\left\{#1\right\}}
\newcommand{\mto}[3]{#1\colon #2 \to #3}
\newcommand{\chu}[1]{(#1,\leqslant_{#1})}
\newcommand{\ls}[1]{\leqslant_{#1}}
\newcommand{\gs}[1]{\geqslant_{#1}}
\newcommand{\uc}[1]{#1^{\Delta}}
\newcommand{\lc}[1]{{#1^{\nabla}}}
\newcommand{\mcX}{\mathcal{X}}

\newcommand{\mfB}{\mathfrak{B}}
\newcommand{\cf}{\textbf{C}_f}
\newcommand{\bb}{\textbf{B}}
\newcommand{\sa}{\textbf{S}_A}
\newcommand{\ZZ}{\mathbb{Z}}

\newcommand{\RR}{\mathbb{R}}

\title{\bf Изотонное продолжение отображений}
\author{\bf А.\,А. Довгошей}
\date{}

\begin{document}
\begin{center}
\textbf{\Large Isotone extension of mappings}

\bigskip

{\bf \large Oleksiy Dovgoshey}
\end{center}
{\bf Abstract.} A set of necessary and sufficient conditions under which an isotone mapping from a subset of a poset $X$  to a poset $Y$ has an extension to an isotone mapping from $X$ to $Y$ are found.

\bigskip

\begin{center}
\textbf{\Large Изотонное продолжение отображений}

\bigskip

{\bf \large А.\,А. Довгошей}
\end{center}

\begin{abstract}
Найден ряд необходимых и достаточных условий того, что изотонное отображение, заданное на подмножестве частично упорядоченного множества $X$ и принимающее значения в частично упорядоченном множестве $Y$, продолжается до изотонного отображения из $X$ в $Y$.
 \\\\
{\bf Ключевые слова и фразы:} изотонное отображение, полная решетка, линейно упорядоченное множество, обобщенная решетка, универсальное частично упорядоченное множество.
\\\\
УДК 512.562, 512.568
\end{abstract}

\section{Введение}
Пусть $\chu{X}$ и $\chu{Y}$~--- частично упорядоченные (ч.\,у.) множества. Отображение $\mto{f}{X}{Y}$ называется \emph{изотонным}, если импликация
\[(x\ls{X}y)\Rightarrow (f(x)\ls{Y} f(y))\]
выполняется для всех $x,y \in X$. В частности, если $f$~--- изотонная биекция,  и $f^{-1}$ изотонно, то $f$~--- \emph{изоморфизм} между $\chu{X}$ и $\chu{Y}$.

Пусть $A\subseteq X$ и $\mto{g}{A}{Y}$ изотонно как отображение ч.\,у. множества $A$ с порядком $\ls{A}$, индуцированным из $X$. Отображение $\mto{\Psi}{X}{Y}$ назовем \emph{изотонным продолжением} отображения $g$, если $\Psi$~---  изотонно и $\Psi(x)=g(x)$ для любого $x\in A$.

Задача построения продолжения отображения $\mto{f}{A}{Y}$, $A\subseteq X$, до изотонного отображения $\mto{g}{X}{Y}$ обычно рассматривается при дополнительных ограничениях на $f$, $g$, $X$, $A$ и $Y$. Например, в~\cite{BNS} изотонное продолжение строится в том случае, когда $X$ и $Y$  замкнутые конусы в топологических векторных пространствах, $A$~--- внутренность конуса $X$, a от $g$ требуется непрерывность или полунепрерывность.

Другая типичная ситуация~--- это продолжение изотонного непрерывного отображения, на упорядоченном или предупорядоченном топологическом пространстве  до изотонного непрерывного отображения на компактификации этого пространства (см., например, ~\cite{Mi} и приведенные там ссылки). Как отмечается в~\cite{Mi}, такого рода исследования мотивированы, в частности, попытками перенести причинно-следственные  отношения (casual relations) на идеальные границы лоренцевых многообразий.

Упомянем также, возникающую в численных методах, задачу монотонный интерполяции монотонных данных (см., например,~\cite{Vol}) и фундаментальную в теории вероятности теорему  о продолжении меры с булевой алгебры  на порожденную ею $\sigma$-алгебру (см., например, ~\cite[\S1.5]{Nev}). Стоит отметить, что в~\cite{Vol} изотонное продолжение исходных изотонных данных осуществляется с помощью  интерполяции кубическими сплайнами, а при продолжении  меры естественно требуется её субаддитивность.

В настоящей работе задача построения изотонного продолжения отображения $\mto{f}{A}{Y}$, $A\subseteq X$, до отображения из $X$ в $Y$ исследуется в чистом виде, без предположения о том, что заданы алгебраические или топологические структуры отличные от частичного порядка.

Доказано, что эта задача разрешима для
\begin{itemize}
  \item любых $X,A \subseteq X$  и всех изотонных $\mto{f}{A}{Y}$ тогда и только тогда, когда $Y$~---  полная решетка (теорема~\ref{t4});
  \item любых $Y$, всех $X$, включающих данное $A$,  и всех изотонных $f$, заданных на $A$ тогда и только тогда, когда $A$~---  полная решетка (теорема ~\ref{t5});
  \item любых $A, Y$ и всех изотонных $\mto{f}{A}{Y}$ тогда и только тогда, когда $X$ изоморфно кардинально неразложимой сумме подмножеств множества $\ZZ$ (теорема~\ref{t10});
  \item $A$, имеющих мощность меньшую данного кардинального числа $\alpha$, любых $X$ и всех изотонных $\mto{f}{A}{Y}$ тогда и только тогда, когда $Y$~--- $\alpha$-квазирешетка (определение~\ref{d4.2*} и теорема~\ref{t4.4*}).
\end{itemize}

Показано, что для любых $X$ и, обладающих наименьшим и наибольшим элементами, $A\subseteq X$ каждое изотонное $\mto{f}{A}{Y}$ продолжается до изотонного отображения, заданного на $X$ и сохраняющего экстремальные значения тогда и только тогда, когда $Y$~--- полная локальная решетка (теорема~\ref{t5.3} и определение~\ref{d5.1}).

\section  {Изотонные продолжения и полные \mbox{решетки}}
Напомним, что ч.\,у. множество $\chu{Y}$ называется \emph{полной решеткой}, если всякое его непустое подмножество $A$ имеет как точную верхнюю так и точную нижнюю грани.  Если существование точных верхней и нижней граней требуется только для конечных $A$, то, по определению, $Y$ является \emph{решеткой}.

Пусть $\chu{Y}$~--- полная решетка, $Y\neq \varnothing$. Наибольший и наименьший элементы из $Y$ будем обозначать  через $1_Y$ и $0_Y$ соответственно,
\[ 1_Y=\sup\nolimits_Y Y \ \text{  и  } \  0_Y= \inf\nolimits_Y Y. \]
Если $A$~--- пустое подмножество $Y$, то
\begin{equation}\label{e1}
\inf\nolimits_Y A=1_Y \ \text{ и } \  \sup\nolimits_Y A = 0_Y.
\end{equation}

Пусть $\chu{Y}$~--- ч.\,у. множество, $A\subseteq Y$. \emph{Верхний конус} множества $A$~--- это подмножество $\uc{A}$ множества $Y$ такое, что
\[y \in \uc{A} \Leftrightarrow a\ls{Y} y \]
для всех $a \in A$ (см., например,~\cite[стр. 11]{Sk}). Двойственным образом определяем \emph{нижний конус} $\lc{A}$. Если $A$ является одноточечным, $A=\set{a}$, то по определению
\[\uc{a} =\uc{A} \ \text{ и } \ \lc{a} =\lc{A}.\]
Элементы $\uc{A}$ называются \emph{мажорантами} множества $A$,  а элементы $\lc{A}$~--- \emph{минорантами} $A$.

\begin{remark}\label{r1}
Хорошо известно, что всякое ч.\,у. множество $\chu{X}$ изоморфно вкладывается в булеан этого множества $\mathfrak{B}(X)$, упорядоченный отношением $\subseteq$,  с помощью отображения
\[ \mto{\nabla_X}{X}{\mathfrak{B}(X)}, \quad  \nabla_X(a)=\lc{a}, \quad a\in X \]
(см., например,~\cite[стр.  41, 42]{Sk}).
\end{remark}

Таким образом, справедлива

\begin{lemma}\label{l2}
Для любого ч.\,у. множества $\chu{Y}$ существует полная решетка $\chu{M}$ такая, что $Y$ изоморфно некоторому множеству $A\subseteq M$.
\end{lemma}

\begin{lemma}\label{l3}
Пусть $\chu{X}$~--- полная решетка и $A\subseteq X$. Если тождественное отображение $\mto{\id}{A}{A}$ имеет изотонное продолжение $\mto{g}{X}{A}$, то $A$~--- полная решетка относительно порядка индуцированного из $X$.
\end{lemma}

\begin{proof}
Пусть $\mto{\id}{A}{A}$ имеет изотонное продолжение $\mto{g}{X}{A}$. Пусть $\varnothing \neq B \subseteq A$. Достаточно установить, что существуют $\sup_A{B}$ и $\inf_A{B}$. Так как $X$~--- полная решетка, то существует $\bar{b}:=\sup_X{B}$. Покажем, что $g(\bar{b})$ представляет собой $\sup_A{B}$. Действительно, так как $x\ls{X}\bar{b}$ выполнено для любого $x\in B$, а $g$~--- изотонное продолжение отображения $\mto{\id}{A}{A}$, то
\[ x=g(x)\ls{A} g(\bar{b}). \]
Следовательно, $g(\bar{b})$~--- мажоранта $B$ в $\chu{A}$. Пусть теперь $t$~---  произвольная мажоранта $B$ в $\chu{A}$. Тогда $t$~--- мажоранта $B$ в $\chu{X}$. Так как $\bar{b}=\sup_X{B}$, то $\bar{b}\ls{X}t$. Следовательно,
\[ g(\bar{b})\leqslant g(t)=\id (t)=t. \]
Теперь из определения точной верхней грани следует, что $g(\bar{b})=\sup_A{B}$. Двойственным образом устанавливается наличие точной нижней грани множества $B$ в $\chu{A}$.
\end{proof}

\begin{definition}\label{d1}
Пусть $X$ и $Y$~--- ч.\,у. множества, $A\subseteq X$, $\mto{f}{A}{Y}$~--- изотонное отображение  и пусть $\mto{\Psi}{X}{Y}$~--- изотонное продолжение отображения $f$. Будем говорить, что $\Psi$~--- \emph{верхнее (нижнее) изотонное продолжение} $f$, если для любого изотонного продолжения $\mto{\mcX}{X}{Y}$ отображения $f$ и всех $x\in X$ имеет место \[ \mcX(x)\ls{Y}\Psi(x)\quad (\Psi(x)\ls{Y} \mcX(x)). \]
\end{definition}

Следующая теорема частично анонсирована в~\cite{DP}, а использованная при её доказательстве конструкция  верхнего изотонного продолжения позаимствована из~\cite{DPK}.

\begin{theorem}\label{t4}
Пусть $\chu{Y}$~--- непустое ч.\,у. множество. Следующие утверждения эквивалентны.
\begin{itemize}
  \item [(i)] $\chu{Y}$~--- полная решетка.
  \item [(ii)] Для любого ч.\,у. множества $\chu{X}$ и любого $A\subseteq X$ каждое изотонное отображение $\mto{f}{A}{Y}$ имеет нижнее изотонное продолжение.
  \item [(iii)] Для любого ч.\,у. множества $\chu{X}$ и любого  $A\subseteq X$ каждое изотонное отображение $\mto{f}{A}{Y}$ имеет верхнее изотонное продолжение.
  \item [(iv)] Для любого ч.\,у. множества $\chu{X}$ и любого $A\subseteq X$ каждое изотонное отображение $\mto{f}{A}{Y}$ имеет изотонное продолжение.
\end{itemize}
\end{theorem}
\begin{proof}
\textbf{(i)}$\Rightarrow$\textbf{(ii)}. Пусть $\chu{Y}$~--- полная решетка, $\chu{X}$~--- ч.\,у. множество, $A\subseteq X$ и $\mto{f}{A}{X}$~--- изотонное отображение. Для любого $x\in X$, положим
\begin{equation}\label{e2}
f_*(x):=\sup\nolimits_Y \{f(t)\colon t\in A\cap \lc{x}\}.
\end{equation}
Заметим, что в соответствии с~(\ref{e1})
\begin{equation}\label{e3}
f_*(x)=0_Y,
\end{equation}
если $A\cap \lc{x}=\varnothing$. Проверим, что отображение \[X\ni x \mapsto f_*(x)\in Y\]
является изотонным. Действительно, $x_1\ls{X}x_2$, то $\lc{x_1}\subseteq \lc{x_2}$, что в силу~(\ref{e2}) дает $f_*(x_1)\ls{Y}f_*(x_2)$. Убедимся в том, что $f_*$ продолжает $f$. Пусть $x\in A$, тогда $t\ls{A}x$ для любого $t\in A\cap \lc{x}$. Отсюда в силу изотонности $f$, находим $f(t)\ls{Y}f(x)$ для всех $t\in A\cap \lc{x}$. Теперь, используя~(\ref{e2}), получаем
\[f_*(x)\ls{Y} f(x).\]
Для доказательства того, что $f_*(x)=f(x)$ осталось установить, что $f(x)\ls{Y} f_*(x)$. Последнее легко следует из~(\ref{e2}), так как
\[(x\in A)\Rightarrow (x\in A\cap\lc{x}),\]
а значит
\[\sup\nolimits_Y\{f(t)\colon t\in A\cap\lc{x}\}\gs{Y}f(x).\]
Таким образом, $f_*$~--- изотонное продолжение $f$. Проверим, что $f_*$~--- нижнее изотонное продолжение $f$. Пусть $\mto{\mcX}{X}{Y}$~--- произвольное продолжение $f$. Нужно установить, что
\begin{equation}\label{e4}
f_*(x)\ls{Y} \mcX(x)
\end{equation}
для любого $x\in X$. Если $A\cap\lc{x}=\varnothing$, то~(\ref{e4}) следует из~(\ref{e3}). Пусть $A\cap\lc{x}\neq\varnothing$. Рассмотрим произвольное $t\in A\cap\lc{x}$. Так как $\mcX$~--- изотонное отображение, то $\mcX(t)\ls{Y}\mcX(x)$. А так как $\mcX$ продолжает $f$  и $t\in X$, то имеем
\[f(t)\ls{Y} \mcX(x)\]
для любого $t\in A\cap\lc{x}$. Отсюда и из ~(\ref{e2}) следует~(\ref{e4}).

\textbf{(i)}$\Rightarrow$\textbf{(iii)}. Проверку этой импликации можно осуществить рассуждениями двойственными рассуждениям приведенным выше. Заметим в частности, что верхнее изотонное продолжение отображения $\mto{f}{A}{Y}$ можно задать формулой
\begin{equation}\label{e5}
f^*(x):=\inf\nolimits_Y{\{f(t)\colon t \in A\cap\uc{x} \}},
\end{equation}
если $\chu{Y}$~--- полная решетка.

Импликации (iii)$\Rightarrow$(iv) и (ii)$\Rightarrow$(iv) очевидны.

\textbf{(iv)$\Rightarrow$(i).} Пусть выполнено (iv). Докажем, что $\chu{Y}$~--- полная решетка. По лемме~\ref{l2} найдется полная решетка $\chu{X}$ и множество $A\subseteq X$ такое, что ч.\,у. множество $\chu{A}$ (с порядком индуцированным из $\chu{X}$) изоморфно $\chu{Y}$. Не уменьшая общности, можно считать $A=Y$. В силу (iv) отображение $\mto{\id}{Y}{Y}$ может быть продолжено до изотонного отображения $\mto{g}{X}{Y}$. По лемме~\ref{l3}, ч.\,у. множество $\chu{Y}$ является полной решеткой.
\end{proof}

Пусть $\chu{Y}$~--- полная решетка, $\chu{X}$~--- ч.\,у. множество, $A\subseteq X$ и $\mto{f}{A}{Y}$~--- изотонное отображение. Обозначим через $\cf$ множество всех изотонных продолжений $\mto{g}{X}{Y}$ отображения $f$. Тогда для нижнего и верхнего изотонных продолжений $f$ вместо формул~(\ref{e2}) и~(\ref{e5}) можно использовать соотношения
\begin{equation}\label{e9}
f_*(x)=\inf\nolimits_Y\{g(x)\colon g\in \cf\}, \qquad f^*(x)=\sup\nolimits_Y\{g(x)\colon g\in \cf\}.
\end{equation}
Действительно, пусть $f_*$ определена первой из формул~(\ref{e9}). Если $x_1,x_2 \in X$ и $x_1\ls{X}x_2$, а $g$--- произвольный элемент $\cf$, то
\[ f_*(x_1)\ls{Y}g(x_1)\ls{Y}g(x_2). \]
Взяв infinum по всем $g\in \cf$, из неравенства $f_*(x_1)\ls{Y}g(x_2)$  получим $f_*(x_1)\ls{Y} f_*(x_2)$, т.\,е. $\mto{f_*}{X}{Y}$~--- изотонное отображение. Равенство $f_*(x)=f(x)$ для $x\in A$ очевидно. Таким образом, $f_*(x)$--- нижнее изотонное продолжение $f$. Аналогично доказываем, что $f^*$~--- верхнее изотонное отображение.

Если на множестве $\cf$ задать частичный порядок как
\[ (g\ls{\cf} \Psi) \Leftrightarrow (g(y)\ls{Y} \Psi(y) \text{ для всех  } y\in Y), \]
то легко получить
\begin{corollary}\label{c2.11*}
Пусть $\chu{Y}$~--- полная решетка. Тогда \linebreak[4] $\chu{\cf}$~--- полная решетка для любого $\chu{X}$, любого $A\subseteq X$ и любого изотонного $\mto{f}{A}{Y}$, а отображения $f_*$ и $f^*$, определенные формулами~(\ref{e2}) и ~(\ref{e5}), являются нулём и, соответственно,  единицей решетки $\cf$,
\begin{equation}\label{e2.11**}
f_*= 0_{\cf}\ \text{ и } \ f^*= 1_{\cf}.
\end{equation}
\end{corollary}
\begin{proof}
Если $\bb$~--- произвольное непустое подмножество множества $\cf$, то для доказательства существования $\inf\nolimits_{\cf}\bb$ можно повторить рассуждения, приведенные при проверке того, что формулы~(\ref{e2}) и~(\ref{e9}) определяют одну и ту же функцию $f_*$, заменив в них $\cf$ на $\bb$. Существование $\sup\nolimits_{\cf}\bb$ доказывается двойственным образом, a равенства~(\ref{e2.11**}) фактически уже доказаны.
\end{proof}

Если $\chu{Y}$~--- решетка и $\varnothing \neq B \subseteq Y$, то подмножество $L$ множества $Y$, состоящее из элементов вида $\sup\nolimits_Y C$ и $\inf\nolimits_Y C$, где $C$~--- конечные непустые подмножества множества $B$ является подрешеткой решетки $Y$. Говорят, что $L$ \emph{порождена} множеством $B$.

Дадим ещё одно следствие теоремы~\ref{t4}.

\begin{corollary}\label{c6}
Пусть $\chu{Y}$~--- решетка. Тогда для любого ч.\,у. множества $\chu{X}$ и любого конечного $A\subseteq X$ любое изотонное отображение $\mto{f}{A}{Y}$ имеет изотонное продолжение.
\end{corollary}

\begin{proof}
При $A= \varnothing$ это очевидно. Пусть $\varnothing\neq A\subseteq X$, $\mto{f}{A}{Y}$~--- изотонное отображение и $B:=f(A)$. Обозначим через $L$ подрешетку решетки $Y$, порожденную множеством $B$. Тогда число элементов $L$ конечно. Следовательно $L$~--- полная решетка. По теореме~\ref{t4} отображение $A\ni x\mapsto f(x)\in L$ продолжается до изотонного отображения $\mto{g}{X}{L}$. Теперь в качестве искомого продолжения из $X$ в $Y$ можно взять  $X\xrightarrow{g}L\xrightarrow{\operatorname{in}} Y$, где $\operatorname{in}(y)=y$ для всех $y\in L$.
\end{proof}

Критерий разрешимости задачи изотонного продолжения отображений $\mto{f}{A}{Y}$ для конечных $A$ будет дан в четвертом разделе работы (см. теорему~\ref{t4.6*}). Отметим только, что существует ч.\,у. множество $\chu{Y}$ не являющееся решеткой, для которого такая задача разрешима при всех $X$ и всех конечных $A\subseteq X$.

Пусть $\chu{A}$~--- частично упорядоченное множество. Обозначим через $\sa$ класс всех ч.\,у. множеств $\chu{X}$, для которых $A\subseteq X$ и порядок $\ls{A}$ совпадает с порядком, индуцированным из $\chu{X}$.

\begin{theorem}\label{t5}
 Произвольное ч.\,у. множество $\chu{A}$ является полной решеткой тогда и только тогда, когда для всех ч.\,у. множеств $\chu{Y}$, любого изотонного отображения $\mto{f}{A}{Y}$ и любого ч.\,у. множества $\chu{X}\in \sa$ существует $\mto{\Psi}{X}{Y}$, изотонно продолжающее $f$.
\end{theorem}

\begin{proof}
Пусть $\chu{A}$~--- полная решетка, $\chu{Y}$~--- ч.\,у. множество, $\mto{f}{A}{Y}$~---  изотонное отображение и $\chu{X}\in \sa$. По теореме~\ref{t4} тождественное отображение $\mto{\id}{A}{A}$ имеет изотонное продолжение $\mto{g}{X}{A}$. Тогда отображение
\begin{equation}\label{e1*}
X\xrightarrow{g}A\xrightarrow{f}Y
\end{equation}
является изотонным продолжением отображения $f$.

Пусть теперь для любого ч.\,у. множества $\chu{Y}$ и любого \mbox{$\chu{X}$} $\in \sa$ любое изотонное отображение $\mto{f}{A}{Y}$ имеет изотонное продолжение $\mto{g}{X}{Y}$.   По лемме~\ref{l2} в $\sa$ найдется $\chu{X}$, являющееся полной решеткой. Возьмем в качестве $Y$ множество $A$. Пусть  порядки $\ls{Y}$ и $\ls{A}$ совпадают и пусть $f$~--- тождественное отображение из $A$  в $A$. Тогда $f$ имеет изотонное продолжение $\mto{g}{X}{A}$ с $X$, являющимся полной решеткой. Следовательно, по лемме~\ref{l3}, $A$~--- полная решетка.
\end{proof}

\section{Цепи и изотонные продолжения}
Пусть $\chu{X}$~--- ч.\,у. множество. Элементы $a,b\in X$ называются \emph{сравнимыми} если $a\ls{X} b$ или $b\ls{X}a$. Определим на $X$ бинарное отношение $\rho$ по правилу: $(a \rho b)$ тогда и только тогда, когда существует конечная последовательность $a_1,\ldots,a_m$ такая, что $a_1=a$,  $a_m=b$ и при $i=1,\ldots, m-1$ элементы $a_i$ и $a_{i+1}$ сравнимы. Тогда $\rho$~--- эквивалентность на $X$ (см., например, \cite[стр. 15]{Sk}). Смежные классы отношения $\rho$ будем называть \emph{кардинально неразложимыми компонентами} ч.\,у. множества $\chu{X}$.  В частности, если отношение $x\rho y$ имеет место для всех $x,y\in X$, то само $X$~--- кардинально неразложимо.

Напомним, что ч.\,у. множество $X$ называется \emph{линейно упорядоченным} или \emph{цепью}, если сравнимы любые два элемента из $X$.

\begin{lemma}\label{l6}
Пусть $\chu{X}$~--- кардинально неразложимое ч.\,у. множество не являющееся цепью. Тогда найдется $A\subseteq X$ такое, что тождественное отображение $\mto{\id}{A}{A}$ не может быть продолжено до изотонного отображения  из $X$ в $A$.
\end{lemma}

\begin{proof}
Так как $X$ не является цепью, то существуют несравнимые между собой $a,b\in X$. Положим $A=\set{a,b}$ и покажем, что тождественное отображение $\mto{\id}{A}{A}$ не имеет изотонного продолжения на всё $X$.

Допустим, что существует изотонное $\mto{g}{X}{A}$, продолжающее $\id$. В силу кардинальной неразложимости $X$ найдется конечная последовательность $a_1,\ldots , a_m$ такая, что $a_1=a$,  $a_m=b$ и при $i=1,\ldots, m-1$ элементы $a_i$ и $a_{i+1}$ сравнимы. Рассмотрим конечную последовательность $g(a_1),\ldots,g(a_m)$. Так как $g$~--- продолжение $\mto{\id}{A}{A}$ и $a_1=a$, то $g(a_1)=a$.  Из определения множества $A$ следует, что $g(a_2)=a$ или $g(a_2)=b$. Покажем, что $g(a_2)=a$. Действительно т.\,к. $a_1$ сравнимо с $a_2$, то $g(a_1)$ сравнимо с $g(a_2)$. Если $g(a_2)=b$, то $b$ сравнимо с $a$, что противоречит их определению. Следовательно, $g(a_2)=a$. Аналогично доказывается, что
\[ a=g(a_3)=\ldots=g(a_{m-1})=g(a_m). \]
Таким образом, $a=g(a_m)=b$, что противоречит несравнимости $a$ и $b$.
\end{proof}

Пусть $\chu{X}$~--- цепь и $x_1,x_2\in X$. Элемент $x_2$ \emph{непосредственно следует} за $x_1$, если $x_1\ls{X}x_2$ и не существует элемента $x\in X$ лежащего строго между $x_1$ и $x_2$, т.е. такого, что $x_1 <_X x <_X x_2$. Аналогично $x_2$ \emph{непосредственно предшествует} $x_1$, если $x_2\ls{X}x_1$ и в $X$ нет элементов, лежащих строго между $x_2$ и $x_1$.

Напомним, что элемент $x_1$ называется \emph{наименьшим} элементом ч.\,у. множества $X$, если $x_1\ls{X}x$ для любого $x\in X$. Двойственным образом определяется \emph{наибольший} элемент в $X$.

\begin{lemma}\label{l7}
Пусть $\chu{X}$~--- цепь, $x_1\in X$ и $A:=\{x\in X:x<_X x_1\}$. Предположим, что $x_1$ не является наименьшим элементом в $X$, но в $X$ нет элементов, непосредственно предшествующих $x_1$. Тогда тождественное отображение $\mto{\id}{A}{A}$ не может быть продолжено до изотонного отображения из $X$ в $A$.
\end{lemma}

\begin{proof}
Допустим, что  $\mto{g}{X}{A}$~--- изотонное продолжение отображения $\mto{\id}{A}{A}$. Положим
\begin{equation}\label{e3.3*}
x_{-1}:=g(x_1).
\end{equation}
Из определения множества $A$ следует, что $x_{-1} \ls{X} x_1$.  Так как $x_{-1}$ не является непосредственным предшественником $x_1$, то для некоторого $x_0\in X$ имеем
\begin{equation}\label{e12}
x_{-1} <_X x_0 <_X x_1.
\end{equation}
Заметим, что $x_0\in A$. Из~(\ref{e12}), (\ref{e3.3*}) и того, что $g$~--- изотонное продолжение отображения $\mto{\id}{A}{A}$ находим
\[ g(x_{-1})\ls{A}g(x_0)\ls{A}g(x_1), \quad x_{-1} \ls{X} x_0 \ls{X} x_{-1}. \]
Значит $x_{-1}=x_0$, что противоречит~(\ref{e12}).
\end{proof}

Следующая лемма есть утверждение двойственное лемме~\ref{l7}.

\begin{lemma}\label{l8}
Пусть $\chu{X}$~--- цепь, $x_1\in X$ и $A:=\{x\in X\colon x>_{X} x_1\}$. Допустим, что $x_1$ не является наибольшим элементом в $X$, но в $X$ нет элементов, непосредственно следующих за $x_1$. Тогда тождественное отображение $\mto{\id}{A}{A}$ не может быть продолжено до изотонного отображения из $X$ в $A$.
\end{lemma}

Обозначим через $\ZZ$ множество всех целых чисел с обычным отношением порядка
\[\ldots -2\leqslant -1\leqslant 0\leqslant1\leqslant2\ldots.\]
Следующая лемма даёт внутреннюю характеристику порядковых типов  подмножеств множества $\ZZ$.

\begin{lemma}\label{l9}
Пусть $\chu{X}$~--- линейно упорядоченно множество. Множество $\chu{X}$ изоморфно подмножеству множества $\ZZ$  тогда и только тогда, когда для любого $B\subseteq X$ и любого $b_1\in B$ выполняется следующее утверждение.
\begin{itemize}
  \item [(i)] Если $b_1$ не является наименьшим элементом в $\chu{B}$, то существует $b_0\in B$ непосредственно предшествующий $b_1$; а если $b_1$ не является наибольшим элементом в $\chu{B}$, то существует $b_2\in B$, непосредственно следующий за $b_1$.
\end{itemize}
\end{lemma}
\begin{proof}
Если $\chu{X}$ изоморфно какому-то подмножеству $\ZZ$, то $(i)$ очевидно.

Обратно, пусть выполнено $(i)$. Докажем, что существует $A\subseteq \ZZ$  изоморфное $\chu{X}$.  Пусть $(x,y)_X$~--- это множество точек $z\in X$, лежащих между $x$ и $y$, т.\,е.
\[(x,y)_X=\{z\in X\colon x<_X z <_X y\}.\]
Определим на $X$ отношение $\sim$ следующим правилом.
Пусть $x,y\in X$, тогда
\[ (x\sim y)\Leftrightarrow ((x=y)\vee((x,y)_X=\varnothing \ \&  \ x<_{X}y)\vee((y,x)_X=\varnothing \ \& \ y<_X x)). \]
Положим $x \approx y$, если существует конечная последовательность $a_1,\ldots,a_m$ такая, что $a_1=x$, $a_m=y$ и $a_i\sim a_{i+1}$ при $i=1,\ldots,{m-1}$. Таким образом, $\approx$ является транзитивным замыканием отношения $\sim$. Так как $\sim$ является рефлексивным и симметричным, то его транзитивное замыкание $\approx$ есть отношение эквивалентности на $X$  (см., например,~\cite[гл. 1, утверждение 4.25]{Ho}). Если $x\approx y$ выполняется для всех $x,y\in X$, то искомый изоморфизм $X$ на подмножество множества $\ZZ$ можно построить следующим образом. Выберем произвольное $x_0\in X$ и поставим ему в соответствие число $0$, $x_0\mapsto 0$. Если $x_0$ не является наибольшим элементом $\chu{X}$, тогда по $(i)$ существует $x_1\in X$, непосредственно следующий за $x_0$. Поставим элементу $x_1$ в соответствие число $1$, $x_1\mapsto 1$. Аналогично, если $x_0$ не является наименьшим элементом $X$, то найдется $x_{-1}$ непосредственно предшествующий $x_0$.  Положим $x_{-1}\mapsto -1$. Аналогично, если $x_1$ не является наибольшим, а $x_{-1}$ наименьшим, то можно найти $x_{-2}$ и $x_2$ такие, что
\[ (x_{-2},x_2)_X\cap X=\set{x_{-1},x_0,x_1} \]
и положить: $x_{-2}\mapsto -2$ , $x_2\mapsto 2$. Не будем формально описывать <<правила остановки>> этого построения, но отметим, что $\chu{X}$ изоморфно множеству $\set{-k,\ldots,-1,0,1,\ldots,m}$ (при изоморфизме $x_i\mapsto i$) тогда и только тогда, когда $x_{-k}$~--- наименьший элемент $\chu{X}$, a $x_m$~--- наибольший элемент этого множества. Можно убедиться и в том, что любому элементу $x\in X$ будет сопоставлено некоторое $m\in\ZZ$. Это следует из определения отношения $\approx$.

Предположим теперь, что не все пары $x,y\in X$ удовлетворяют отношению $\approx$. Доказательство будет закончено, если показать, что это предположение противоречит $(i)$. Перед тем как осуществить намеченное reductio ad absurdum установим следующее свойство:

(S) <<Если $A_1$ и $A_2$~--- два различных смежных класса, порожденных отношением $\approx$, $a_1 \in A$, $a_2\in A_2$ и $a_1 <_X a_2$, то для любого $y\in A_2$ выполнено $a_1<_X y$>>.

Действительно, пусть существует $y\in A_2$ такое, что
\[ y <_X a_1 <_X a_2, \]
т.\,е. $a_1\in (y,a_2)_X$. Последнее не может выполняться так как $y\approx a_2$, а для любых  $w,u,v \in X$ соотношения
\[ u\approx v,\quad u<_X v, \quad w\in (u,v)_X \]
влекут $u\approx w\approx v$, что легко  доказать индукцией по $m$, где $m$~--- длина самой короткой последовательности $c_1,\ldots,c_m$, для которой $c_1=u$, $c_m=v$ и $c_i\sim c_{i+1}$, $i=1,\ldots, {m-1}$.

Пусть $\set{A_i\colon i\in I}$~--- множество различных смежных классов, порожденных $\approx$. В соответствии с предположением имеем $\abs{I}\geqslant 2$. Для каждого $i\in I$ выберем элемент $a_i\in A_i$. Пусть $A=\set{a_i\colon i\in I}$. Так как $A$~--- цепь, лежащая в $\chu{X}$ и $\abs{I} \geqslant 2$, то, по свойству (i), существуют индексы $i_1$ и $i_2$ такие, что $a_{i_1}$ непосредственно предшествует $a_{i_2}$. Покажем, что в цепи $A_{i_1}$ нет наибольшего элемента. Действительно, пусть $a_{i_1}^*$~--- наибольший элемент в $A_{i_1}$. Так как $a_{i_1}, a_{i_1}^*\in A_{i_1}$, $a_{i_2}\in A_{i_2}$ и $a_{i_1}<_X a_{i_2}$, то по свойству (S) $a_{i_1}^*<_X a_{i_2}$. Отсюда по свойству (i) следует, что существует $b \in X$, непосредственно следующий за $a_{i_1}^*$ в $X$. Тогда $b\sim a_{i_1}^*$, а, следовательно, $b\in A_{i_1}$ и $a_{i_1}^*<_X b$. Это противоречит тому, что $a_{i_1}^*$~--- наибольший элемент цепи $A_{i_1}$.

Рассмотрим теперь цепь $L=A_{i_1}\cup\set{a_{i_2}}$. Так как $a_{i_1}, a_{i_2}\in L$ и $a_{i_1}<_X a_{i_2}$, то $a_{i_2}$ не является наименьшим элементом цепи $L$.  Следовательно, по свойству (i) в $L$ существует $b_{i_1}$, непосредственно предшествующий $a_{i_2}$ в $L$,
\begin{equation}\label{et1}
(b_{i_1}, a_{i_2})_L=\varnothing.
\end{equation}
Так как $b_{i_1}\neq a_{i_2}$, то из определения $L$ следует, что $b_{i_1}\in A_{i_1}$. Так как в $A_{i_1}$ нет наибольшего элемента, то найдется $c_{i_1}\in A_{i_1}$ такой, что
\begin{equation}\label{et2}
b_{i_1}<_X c_{i_1}.
\end{equation}
Используя свойство (S), легко показать, что $a_{i_2}$~--- наибольший элемент цепи $L$. Следовательно,
\begin{equation}\label{et3}
c_{i_1} <_X a_{i_1}.
\end{equation}
Из~(\ref{et2}) и~(\ref{et3})  получим $c_{i_1}\in (b_{i_1}, a_{i_2})_L$, что противоречит~(\ref{et1}).
\end{proof}

\begin{remark}
Условие (i) леммы~\ref{l9} эквивалентно тому, что множество $B$ содержит наибольший элемент, если $B$ ограничено сверху в $X$, и содержит наименьший элемент, если это множество является ограниченным снизу.
\end{remark}

Как показано в \cite{E, H} порядковые типы  подмножеств множества $\ZZ$ играют большую роль при исследовании рассеянных (scattered) линейно упорядоченных множеств. Лемму~\ref{l9} можно вывести из частного случая теоремы 5.37 монографии~\cite{R}, характеризующей рассеянные множества с помощью так называемого $F$-ранга, но данное выше доказательство этой леммы более элементарно.

\begin{theorem}\label{t10}
Пусть $\chu{X}$~--- ч.\,у. множество. Следующие утверждения эквивалентны.
\begin{itemize}
  \item [(i)] Каждая кардинально неразложимая компонента $B$ множества $X$ изоморфна некоторому $C_B\subseteq \ZZ$.
  \item [(ii)] Для любого $A\subseteq X$ и любого ч.\,у. множества  $\chu{Y}$ любое изотонное отображение $\mto{f}{A}{Y}$  может быть продолжено до изотонного отображения $\mto{g}{X}{Y}$.
\end{itemize}
\end{theorem}

\begin{proof}
\textbf{(i)}$\Rightarrow$\textbf{(ii)}.  Пусть выполнено (i) и $\varnothing \neq A \subseteq X$. Для доказательства (ii) достаточно проверить, что тождественное отображение  $\mto{\id}{A}{A}$ имеет изотонное продолжение $\mto{g}{X}{A}$. В этом случае, как было отмечено при доказательстве теоремы~\ref{t5},  отображение $X\xrightarrow{g}A\xrightarrow{f}Y$ есть изотонное продолжение изотонного отображения $\mto{f}{A}{Y}$.

Рассмотрим сначала случай, когда само $\chu{X}$ является кардинально неразложимым. В силу (i) мы можем предположить, что $X\subseteq \ZZ$. Для любого $A\subseteq X$ выполнены следующие альтернативы:
\begin{itemize}
  \item [($i_1$)] в $A$ нет ни наибольшего ни наименьшего элементов;
  \item [($i_2$)] в $A$ есть наибольший элемент, но нет наименьшего;
  \item [($i_3$)] в $A$ есть наименьший элемент, но нет наибольшего;
  \item [($i_4$)] в $A$ есть и наибольший и наименьший элементы.
\end{itemize}
В случаях  ($i_1$) и ($i_2$) искомое изотонное отображение $\mto{g}{X}{A}$ можно задать как
\begin{equation}\label{ew1}
g(x) :=\sup\nolimits_A(\lc{x}\cap A)
\end{equation}
для любого $x\in X$. Так как $A\subseteq X\subseteq \ZZ$, то проверка того, что $g$~--- изотонное продолжение $\mto{\id}{A}{A}$ становится тривиальной. В случае, если выполнено ($i_3$) или ($i_4$), то отображение $g$, заданное формулой ~(\ref{ew1}), будет изотонным продолжением  $\mto{f}{A}{A}$, если как в~(\ref{e1}) считать, что $\sup_A(\varnothing)$~--- представляет собой наименьший элемент множества $A$. Таким образом, (ii) установлено, если $X$~--- кардинально неразложимое множество.

Если $X$ не является кардинально неразложимым, то представим $X$ как кардинальную сумму кардинально неразложимых компонент $X_{\alpha}$, $\alpha\in I$, где $I$~--- некоторое множество индексов  с $\abs{I}\geqslant 2$ (см., например, \cite[стр. 15]{Sk}).

Пусть $x_0$~--- произвольная точка множества $A$  и $A_{\alpha}:=A\cap X_{\alpha}$ при всех $\alpha \in I$. Определим $\mto{g}{X}{A}$ правилом:
\begin{equation}\label{ew2}
g(x)=
\begin{cases}
\sup_{A_{\alpha}}(\lc{x}\cap A_{\alpha}) &\text{при } x\in X_{\alpha} \text{ и } A_{\alpha}\neq \varnothing,\\
x_0 &\text{при } x\in X_{\alpha} \text{ и } A_{\alpha}= \varnothing.
\end{cases}
\end{equation}
Так как $\set{X_{\alpha}\colon \alpha\in I}$~--- разбиение множества $X$, то $g(x)$ определено для любого $x\in X$. В соответствии с доказанным выше сужение $g|_{X_{\alpha}}$ является изотонным продолжением $\id |_{A_{\alpha}}$ при любом $\alpha \in I$. (При $A_{\alpha}=\varnothing$ отображение $\id |_{A_{\alpha}}$ пусто как отображение из пустого множества в $A_{\alpha}$ и значит любое изотонное отображение $X_{\alpha}\to A$  является изотонным продолжением отображения $\id |_{A_{\alpha}}$.)  Осталось заметить, что само отображение $g$ является изотонным. В самом деле, если $x \ls{X} y$, то $x$ сравнимо с $y$, а значит существует такой индекс $\alpha \in I$, для которого $x,y \in X_{\alpha}$. Так как $g_{\alpha}$~--- изотонно, то
\[ (x\ls{X}y)\Rightarrow (x\ls{X_{\alpha}}y) \Rightarrow (g|_{X_{\alpha}}(x)\ls{A}g|_{X_{\alpha}}(y)) \Rightarrow (g(x)\ls{A}g(y)) \]
\textbf{(ii)}$\Rightarrow$\textbf{(i)}. Пусть (ii) выполнено и пусть $X_{\alpha}$~--- произвольная кардинально неразложимая компонента $X$. Тогда для любого $A\subseteq X_{\alpha}$ отображение $\mto{\id}{A}{A}$ продолжается до изотонного отображения на $X$, а значит и на $X_{\alpha}$. Отсюда по лемме~\ref{l6} следует, что $X_{\alpha}$~--- цепь. Используя леммы~\ref{l7} ~\ref{l8} с $X=X_{\alpha}$, убеждаемся в том, что для любого $B\subseteq X_{\alpha}$ и любого $b\in B$ выполнено (i) из леммы~\ref{l9}. Наконец, применяя лемму~\ref{l9}, получим утверждение (i) настоящей теоремы.
\end{proof}

\section{Изотонное продолжение отображений \\ с \mbox{подмножеств} ограниченной мощности}
В настоящем разделе задача изотонного продолжения отображения \linebreak[4] $\mto{f}{A}{Y}$ рассматривается при условии
\begin{equation}\label{e4.1*}
|A|<\alpha,
\end{equation}
где $|A|=\card A$ и $\alpha$~--- фиксированное кардинальное число. В частности при $\alpha=\aleph_0$, где как обычно $\aleph_0$~--- первый бесконечный кардинал, условие~(\ref{e4.1*}) равносильно конечности $A$. При $\alpha>\aleph_0$ задача изотонного продолжения при условии~(\ref{e4.1*}) является, в определенном смысле, менее элементарной. Поэтому начнем со случая $\alpha=\aleph_0$.

В соответствии со следствием~\ref{c6} следует, что задача изотонного  продолжения отображения $\mto{f}{A}{Y}$, $A\subseteq X$ до изотонного отображения $\mto{g}{X}{Y}$ для  конечных $A$ разрешима, если  $Y$~---  решетка (частный случай следствия~\ref{c6}). Таким образом для того, чтобы дать \textbf{критерий} разрешимости такой задачи необходимо подходящее обобщение понятия решетки.

\begin{definition}\label{d4.1}
Будем говорить, что ч.\,у. множество $\chu{Y}$ является \emph{квазирешеткой}, если для любых конечных $A, B\subseteq Y$, удовлетворяющих условиям
\begin{equation}\label{e4.1}
A\subseteq \lc{B} \ \text{ и } \ B\subseteq \uc{A},
\end{equation}
найдется $y^*=y^*(A,B)$ такое, что
\begin{equation}\label{e4.2}
a\ls{Y} y^* \ls{Y} b
\end{equation}
при всех $a\in A$ и $b\in B$.
\end{definition}

\begin{remark}\label{r4.2}
Положив в~(\ref{e4.1}) $A=\varnothing$, получим тривиальные включения $\varnothing \subseteq \lc{B}$,  $B\subseteq Y$. Следовательно, если $\chu{Y}$~--- квазирешетка, то для любого конечного $B\subseteq Y$ существует миноранта. Аналогично проверяется существование мажорант у конечных $B\subseteq Y$.
\end{remark}

\begin{remark}\label{r4.1*}
Неравенство~(\ref{e4.2}) выполняется для всех $a\in A$ и всех $b\in B$ тогда и только тогда, когда оно выполняется для всех максимальных элементов $a$ множества $A$ и всех минимальных элементов $b$ множества $B$. Следовательно, в определении~\ref{d4.1} вместо любых конечных  $A,B \subseteq Y$, удовлетворяющих условию~(\ref{e4.1}), достаточно брать любые конечные антицепи $A,B \subseteq Y$, удовлетворяющие этому условию.
\end{remark}

Очевидно, что любая решетка является квазирешеткой. В обратном направлении верно следующее

\begin{statement}\label{s4.3}
Любая конечная квазирешетка является решеткой.
\end{statement}

\begin{proof}
Пусть $\chu{Y}$~--- конечная квазирешетка и пусть $A \subseteq Y$. Положим $B=\uc{A}$. Как известно, имеет место включение $A\subseteq A^{\Delta \nabla}$ (см., например, \cite[стр. 11]{Sk}). Следовательно,
\[ A\subseteq \lc{B}  \ \text{ и } \ B \subseteq \uc{A}. \]
Так как $\chu{Y}$~--- квазирешетка, то найдется $y^*=y^*(A,B)$, для которого~(\ref{e4.2}) выполнено при всех $a\in A$ и всех $b \in B=\uc{A}$. Значит $y^*=\sup\nolimits_Y A$. Аналогично доказывается существование $\inf\nolimits_Y A$.
\end{proof}

Пример  бесконечной квазирешетки, не являющейся решеткой, будет дан ниже после доказательства теоремы~\ref{t4.6*}.

\begin{lemma}\label{t4.4}
Пусть $\chu{Y}$~--- непустое ч.\,у. множество. Для того, чтобы $\chu{Y}$ было квазирешеткой необходимо и достаточно, чтобы для любого конечного ч.\,у. множества $\chu{X}$  и любого $A\subseteq X$ каждое изотонное отображение $\mto{f}{A}{Y}$ продолжалось до изотонного отображения из $X$ в $Y$.
\end{lemma}
\begin{proof}
\textbf{Необходимость.} Пусть $\chu{Y}$~--- квазирешетка. Покажем, что для любого конечного $\chu{X}$ и любого $A\subseteq X$ любое изотонное отображение $\mto{f}{A}{Y}$ имеет изотонное продолжение. Используем индукцию по $\abs{X}$. Наличие желаемого продолжения очевидно при $\abs{X}=1$. Предположим, что такое продолжение существует при всех $\abs{X}\leqslant m$. Пусть $\abs{X}=m+1$ и $\mto{f}{A}{Y}$ изотонно, $A\subseteq X$. Докажем существование изотонного продолжения $\mto{g}{X}{Y}$  для $f$. Если $A=X$, то все доказано. Если $A\neq X$, то пусть $x_1\in X\setminus A$ и $A_1:=X\setminus \set{x_1}$. По предположению индукции существует изотонное продолжение $\mto{g_1}{A_1}{Y}$ отображения $f$. Будем считать, что $g(x)=g_1(x)$ при всех $x\in A_1$. Осталось определить $g(x_1)$. Пусть $\underline{X}_1:=\lc{x_1}\setminus\set{x_1}$ и $\overline{X}_1:=\uc{x_1}\setminus\set{x_1}$. Положим
\[ \underline{B} := g_1(\underline{X}_1) \ \text{ и } \ \overline{B}:=g_1(\overline{X}_1). \]
Тогда в силу изотонности отображения $g_1$, неравенство $\underline{b}\ls{Y}\overline{b}$ выполняется для всех $\underline{b}\in \underline{B}$ и всех  $\overline{b}\in \overline{B}$. Теперь, используя то, что $\chu{Y}$~--- квазирешетка, можно найти $y^*\in Y$ для которого
\begin{equation}\label{e4.3}
\underline{b}\ls{Y} y^* \ls{Y}\overline{b}
\end{equation}
при всех $\underline{b}\in \underline{B}$, $\overline{b}\in \overline{B}$. Доопределим отображение $g$ в точке $x_1$, положив
\[ g(x_1):=y^*. \]
Очевидно, что $g$~--- продолжение $f$ на $X$. Проверим изотонность этого продолжения. В силу изотонности $g_1$ достаточно установить импликации
\[ ((x<_X x_1)\Rightarrow (g(x)\ls{Y}y^*)) \ \text{и} \ ((x_1<_X x)\Rightarrow (y^*\ls{Y}g(x))), \]
что следует из~(\ref{e4.3}). Действительно, если $x<_X x_1$, то $x\in \underline{X}_1$, a значит  $g(x)\in \underline{B}$. Так как~(\ref{e4.3}) справедливо при всех $\underline{b}\in \underline{B}$ мы имеем $g(x)\ls{Y} y^*$ . Аналогично проверяется и вторая импликация. Таким образом, если $\chu{Y}$~--- квазирешетка, то искомое изотонное продолжение существует.

\textbf{Достаточность.} Пусть для любого конечного $\chu{X}$ и каждого $A\subseteq X$  любое изотонное отображение $\mto{f}{A}{Y}$ продолжается до отображения изотонного на $X$. Докажем, что $\chu{Y}$~--- квазирешетка. Пусть $A$   и $B$~--- два конечных подмножества $Y$, для которых выполняется~(\ref{e4.1}). Достаточно найти $y^*$ такое, что~(\ref{e4.2})  имеет место при всех $a \in A$ и $b\in B$. Пусть $C=A\cup B$ и $\chu{L}$~--- конечная структура, для которой $C\subseteq L$ и такая, что отношения порядка, индуцированные на $C$ из $\chu{Y}$ и из $\chu{L}$  совпадают (см. лемму~\ref{l2}). По теореме~\ref{t4} отображение $\mto{\operatorname{in}}{C}{Y}$, $\operatorname{in}(c)=c$ при всех $c\in C$, продолжается до изотонного отображения $\mto{g}{L}{Y}$. Пусть $\bar{a}=\sup\nolimits_L A$. Тогда
\begin{equation}\label{e4.4}
a\ls{L}\bar{a} \ls{L} b
\end{equation}
при всех $a\in A$, $b\in B$. Теперь из изотонности $g$ и ~(\ref{e4.4}) получаем
\[ a=g(a) \ls{Y} g(\bar{a})\ls{Y}g(b)=b, \]
что дает~(\ref{e4.2}) с $y^*=g(\bar{a})$.
\end{proof}

Обозначим через $FC$ множество всех конечных подмножеств счетно бесконечного множества.

\begin{theorem}\label{t4.6*}
Пусть $\chu{Y}$~--- непустое ч.\,у. множество. Следующие утверждения эквивалентны.
\begin{itemize}
  \item [(i)] $\chu{Y}$~--- квазирешетка.
  \item [(ii)] Для любого конечного ч.\,у. множества $\chu{X}$ и любого $A\subseteq X$ каждое изотонное отображение $\mto{f}{A}{Y}$ изотонно продолжается на $X$.
  \item [(iii)] Для любого ч.\,у. множества $\chu{X}$ и любого конечного $A\subseteq X$ каждое изотонное отображение $\mto{f}{A}{Y}$ изотонно продолжается на $X$.
  \item [(iv)] Для любого конечного подмножества $A$ ч.у. множества $(FC,\subseteq)$ каждое изотонное отображение $\mto{f}{A}{Y}$ изотонно продолжается на $FC$.
\end{itemize}
\end{theorem}

\begin{proof}
Эквивалентность (i)$\Leftrightarrow$(ii) установлена в лемме~\ref{t4.4}. Импликации (iii)$\Rightarrow$(ii) и (iii)$\Rightarrow$(iv) очевидны.

Проверим справедливость (ii)$\Rightarrow$(iii). Пусть выполнено (ii). Рассмотрим произвольное $\chu{X}$, конечное $A\subseteq X$ и изотонное $\mto{f}{A}{X}$. Обозначим через $LX$ полную решетку, для которой $X\subseteq LX$ и такую, что порядок, индуцированный из $LX$ на $X$, совпадает с $\ls{X}$. Пусть $LA$~--- подрешетка решетки $LX$, порожденная множеством $A$. Так как $LA$~--- конечное ч.\,у. множество, то по лемме~\ref{t4.4} существует изотонное отображение $\mto{g}{LA}{Y}$, продолжающее отображение $\mto{f}{A}{Y}$. Так как любая конечная решетка является полной, то по теореме~\ref{t4} существует изотонное продолжение $\mto{\Psi}{LX}{LA}$ вложения $\mto{\operatorname{in}_1}{LA}{LA}$, где $\operatorname{in}_1(a)=a$ для любого $a\in LA$. Пусть $\mto{\operatorname{in}_2}{X}{LX}$~--- вложение $X$  в $LX$, $\operatorname{in}_2(x)=x$ для всех $x\in X$. Легко проверить, что отображение
\[ X \xrightarrow{\operatorname{in}_2} LX \xrightarrow{\Psi} LA \xrightarrow{g} Y\]
есть одно из искомых изотонных продолжений отображения $f$.

Для завершения доказательства остается установить, что имеет место (iv)$\Rightarrow$(ii). Сделаем это. Пусть выполнено (iv) и пусть $\chu{X}$~--- конечное ч.у. множество, $A\subseteq X$ и $\mto{f}{A}{Y}$~--- изотонное отображение. Пусть $X_F$~--- подмножество $FC$ изоморфное множеству $\chu{X}$.  Тогда существуют $A_F\subseteq X_F$ изоморфное $A$ и изотонное отображение $\mto{f_F}{A_F}{Y}$ такие, что коммутативна диаграмма
\[ \begin{diagram}
\node{A} \arrow[1]{s,l}{\operatorname{is}}
\arrow{se,t}{f}\\
\node{A_F}
\arrow{e,b}{f_F}
\node{Y,}
\end{diagram}\]
где $\operatorname{is}$~--- изоморфизм между $A$ и $A_F$. Достаточно найти изометричное продолжение $f_F$ на $X_F$. В силу (iv) существует $\mto{g}{FC}{Y}$, изотонно продолжающее $f_F$. Взяв сужение $g$ на $X_F$ получаем искомое.
\end{proof}

\begin{example}\label{ex7}
Пусть $\RR$~---множество всех действительных чисел, наделенных обычным порядком $\leqslant$. Положим
\[R=(\RR\setminus\set{0})\cup\set{0_1,0_2},\]
где $0_1$ и $0_2$ выбраны так, что  $(\RR\setminus\set{0})\cap\set{0_1,0_2}=\varnothing$. Определим на $R$  порядок следующим образом:

если $x,y\in \RR\setminus\set{0}$, то $(x\ls{R}y)\Leftrightarrow(x\leqslant y)$;

если $x\in \RR\setminus\set{0}$, $y\in \set{0_1,0_2}$, то
\[ (x\ls{R}y)\Leftrightarrow(x\leqslant 0), \ \text{ и } \ (y\ls{R}x)\Leftrightarrow(0\leqslant x); \]

если $x,y\in\set{0_1,0_2}$, то $(x\ls{R}y)\Leftrightarrow (x=y)$.

Ч.\,у. множество $\chu{R}$ не является решеткой, так как множество $\set{0_1,0_2}$ не имеет верхней грани в $R$. Непосредственно проверяется, что $\chu{R}$~--- квазирешетка. Следовательно,  по теореме \ref{t4.6*}, для любого конечного $\chu{X}$ и любого $A\subseteq X$ каждое изотонное отображение $\mto{f}{A}{R}$ продолжается до изотонного отображения на $X$.
\end{example}

Из теорем~\ref{t4.6*} и~\ref{t4} вытекает следующее

\begin{corollary}\label{c4.6}
Пусть $\chu{Y}$~--- конечно. Предположим задача изотонного продолжения изотонных отображений $A\to Y$ разрешима для всех конечных ч.\,у. множеств $\chu{X}$ и всех $A\subseteq X$. Тогда  эта задача разрешима для всех $\chu{X}$ и всех $A\subseteq X$.
\end{corollary}

\begin{proof}
По теореме~\ref{t4.6*}, $\chu{Y}$~--- квазирешетка. В силу утверждения~\ref{s4.3} эта квазирешетка является решеткой.  Так как эта решетка конечна, то она полна. Остается применить теорему~\ref{t4}.
\end{proof}

Перейдем теперь к проблеме изотонного продолжения отображений, заданных на $A$, удовлетворяющих~(\ref{e4.1}) с произвольным $\alpha\geqslant \aleph_0$.

\begin{definition}\label{d4.2*}
Пусть $\alpha$~--- бесконечный кардинал. Ч.у. множество $\chu{Y}$ является \emph{$\alpha$-квазирешеткой}, если для любых $A,B \subseteq Y$, удовлетворяющих включениям ~(\ref{e4.1}) и соотношениям $|A|<\alpha$, $|B|<\alpha$ найдется $y^*=y^*(A,B)\in Y$ такое, что ~(\ref{e4.2}) выполнено при всех $a\in A$ и $b\in B$.
\end{definition}

В такой терминологии квазирешетки~---  это в точности $\aleph_0$- квазирешетки. Очевидно, что любая $\alpha$-квазирешетка будет $\beta$-квазирешеткой при $\beta\leqslant\alpha$. В частности все $\alpha$-квазирешетки будут квазирешетками.

Квазирешетка из примера~(\ref{ex7}) является $\aleph_0$-квазирешеткой, но не является $\mathfrak{c}$-квазирешеткой (здесь, как обычно, $\mathfrak{c}$ это мощность континуума).

Следующее определение является <<двойником>> определения 2.1 из главы 5 монографии~\cite{Ha}.

\begin{definition}\label{d4.3*}
Пусть $\beta$~--- бесконечный кардинал. Ч.у.  множество $\chu{P}$ назовем $\beta$-универсальным, если для любого ч.у. множества \mbox{$\chu{X}$} с $|X|<\beta$ существует $T\subseteq P$ такое, что $\chu{T}$ изоморфно $\chu{X}$.
\end{definition}

При доказательстве теоремы~\ref{t4.4*}, которая является основным результатом настоящего раздела, будет предполагаться, что:
\begin{itemize}
  \item каждый бесконечный кардинал $\beta$ отождествлен с наименьшим ординалом, имеющим мощность $\beta$, т.е. $\beta=\aleph_{\alpha}$, где $\alpha$~--- соответствующее ординальное число;
  \item $2^{\aleph_{\alpha}}=\aleph_{\alpha+1}$ при любом ординале $\alpha$, т.е. справедлива обобщенная континуум-гипотеза (GCH);
  \item ординальные числа~--- это транзитивные множества вполне упорядоченные отношением $\in$. В частности $\beta=\aleph_{\alpha}$ есть вполне упорядоченное множество ординалов мощности строго меньшей чем $\beta$ (см., например,~\cite[Гл. 2]{J}).
\end{itemize}

Кроме того, при построении изотонных продолжений отображений $\mto{f}{A}{Y}$  будем использовать следующую конструкцию. Пусть $I$~--- цепь и пусть каждому $i\in I$ сопоставлено множество $A_i\supseteq A$ и отображение  $\mto{g_i}{A_i}{Y}$ такие, что $g_i|_A=f$ и $A_i\subseteq A_j$ и $g_j|_{A_i}=g_i$ при $i\ls{I}j$. Тогда на множестве $X:=\bigcup\limits_{i\in I}A_i$ естественным образом определено единственное отображение $\mto{g}{X}{Y}$, для которого $g|_{A_i}=g_i$ при всех $i \in I$. Будем обозначать это отображение $g$ символом $\Lim\limits_{i\in I}g_i$. Легко видеть, что если $X$ наделено порядком $\ls{X}$, а все $A_i$ порядками $\ls{A_i}$ такими, что
\[ \ls{A_i} = (A_i\times A_i)\cap \ls{X} \]
и каждое $\mto{g_i}{A_i}{Y}$~--- изотонно, то $\Lim\limits_{i\in I}g_i$~--- изотонное продолжение $f$ на $X$.

\begin{theorem}\label{t4.4*}
Пусть $\chu{Y}$~--- непустое ч.у. множество и $\alpha$~--- бесконечный кардинал и $\chu{P}$~--- $\alpha$-универсально. Следующие утверждения являются эквивалентными.
\begin{itemize}
  \item [(i)] $\chu{Y}$~--- $\alpha$-квазирешетка.
  \item [(ii)] Для любого ч.у. множества $\chu{X}$ с $|X|<\alpha$ и любого $A\subseteq X$ каждое изотонное отображение $\mto{f}{A}{Y}$ изотонно продолжается на $X$.
  \item [(iii)] Для любого ч.у. множества $\chu{X}$  и любого $A\subseteq X$  с $|A|<\alpha$ любое изотонное отображение $\mto{f}{A}{Y}$ изотонно продолжается на $X$.
  \item [(iv)] Для любого $A\subseteq P$ с $|A|<\alpha$ каждое изотонное отображение $\mto{f}{A}{Y}$ изотонно продолжается на $P$.
\end{itemize}
\end{theorem}

\begin{proof}
\textbf{(i)$\Rightarrow$ (ii).} Пусть $\chu{Y}$~--- $\alpha$-квазирешетка, $\chu{X}$~--- ч.у. множество с $|X|<\alpha$, $A\subseteq X$, $X\setminus A\neq\varnothing$ и $\mto{f}{A}{Y}$~--- изотонно. Используя трансфинитную индукцию, докажем, что $f$ изотонно продолжается на $X$.  Зададим  на $X\setminus A$ порядок $\preccurlyeq$, превращающий $X\setminus A$ во волне упорядоченное множество. Пусть $x_i$~---  такой элемент  из $(X\setminus A,\preccurlyeq)$, что для любого $x_j\prec x_i$ отображение $f$ изотонно продолжается до отображения $\mto{g_{x_j}}{A\cup A_j}{Y}$, где
\[ A_j:=\{a\in X\setminus A\colon x\preccurlyeq x_j \}, \]
так, что $g_{x_j}|_{A_k}=g_{x_k}$, если $x_k\preccurlyeq x_j$. Множество таких $x_i$ не пусто, в частности ему принадлежит наименьший элемент из $(X\setminus A,\preccurlyeq)$. Положим
\begin{equation}\label{s1}
g_{x_i}^{\circ}:=\Lim\limits_{X_j\in A_i^{\circ}}g_{x_j},
\end{equation}
где $A_i^{\circ}:=\{x\in X\setminus A\colon x\prec x_i \}$. Тогда $\mto{g_{x_i}^{\circ}}{A\cup A_{i}^{\circ}}{Y}$~--- изотонное продолжение $f$. Продолжим $g_{x_i}^{\circ}$ до изотонного отображения $\mto{g_{x_i}}{A\cup A_i}{Y}$. Пусть
\[ \underline{x}:=\lc{x_i}\cap (A\cup A_i^{\circ}) \ \text{ и } \   \overline{x}:=\uc{x_i}\cap (A\cup A_i^{\circ}),  \]
где $\lc{x_i}$  и $\uc{x_i}$~---  есть нижний и, соответственно, верхний конусы элемента $x_i$ в $\chu{X}$. Очевидно, $\underline{x}\ls{X}x_i\ls{X}\overline{x}$ для любых $\overline{x}\in \overline{X}_i$  и $\underline{x}\in \underline{X}_i$. Положим $\underline{B}:=g_{x_i}^{\circ}(\underline{X}_i)$,  $\overline{B}:=g_{x_i}^{\circ}(\overline{X}_i)$. Тогда $\underline{b}\ls{Y}\overline{b}$  для всех $\underline{b} \in \underline{B}$ и $\overline{b}\in \overline{B}$. Кроме того очевидно
\[ |\underline{B}\cup \overline{B}|\leqslant |\underline{X}_i\cup \overline{X}_i|\leqslant |X|<\alpha. \]
Так как $\chu{Y}$~--- $\alpha$-квазирешетка, то существует $y_*\in Y$  такое, что
\[ \underline{b}\ls{Y}y_* \ls{Y} \overline{b} \]
при всех $\underline{b} \in  \underline{B}$  и  $\overline{b} \in  \overline{B}$. Положим
\[ g_{x_i}(x):=
\begin{cases}
g_{x_i}^{\circ}(x) &\text{ при } \ A\cup A_i^{\circ}  \\
y^* &\text{ при }  \ x=x^*.
\end{cases}
 \]
Непосредственно проверяется, что $g_{x_i}$~--- изотонно и
\begin{equation}\label{e4.2*}
g_{x_k} = g_{x_i}|_{A\cup A_{x_k}} \ \text{ при } \ x_k \preccurlyeq x_i.
\end{equation}
Следовательно, в соответствии с принципом трансфинитной индукции, для любого $x_i\in X\setminus A_j$ существует $\mto{g_{x_i}}{A\cup A_i}{Y}$, изотонно продолжающее $f$ и удовлетворяющее~(\ref{e4.2*}) при $x_k\preccurlyeq x_i$. Осталось положить
\begin{equation}\label{s2}
g:=\Lim\limits_{x_i\in X\setminus A}g_{x_i}
\end{equation}
и мы получаем искомое изотонное продолжение $f$ на $X$.

\textbf{(ii)$\Rightarrow$(iii).} Пусть выполнено (ii), $\chu{X}$~--- ч.у. множество, $A\subseteq X$, $|A|<\alpha$  и $\mto{f}{A}{Y}$~--- изотонно. Докажем, что $f$  изотонно продолжается на $X$. Если $|A|=|X|$, то все доказано. Предположим вначале, что $|X|$~--- первый кардинал строго больший чем $|A|$. Пусть $\beta$~--- наименьший ординал, имеющий мощность $|X|$.  Тогда на $X\setminus A$  можно задать вполне-упорядочение $\preccurlyeq$ такое, что $(X\setminus A,\preccurlyeq)$ изоморфно $\beta$. Как и при доказательстве импликации (i)$\Rightarrow$(ii) рассмотрим элемент $x_i\in X\subseteq A$ такой, что $f$ продолжается на $A\cup A_j$ при $j\prec i$ и продолжим $f$ на $A^{\circ}_i$ как в~(\ref{s1}). (Используются обозначения из первой части доказательства настоящей теоремы.) Вспоминая соглашения, приведенные перед формулировкой теоремы, убеждаемся в том, что $|A^{\circ}_i|<\alpha$. В этом случае $|A\cup A^{\circ}_i|= |A|+|A^{\circ}_i|<\alpha$. Следовательно, в силу (i), найдется $\mto{g_{x_i}}{A\cup A_i}{Y}$, изотонно продолжающее $f$. Искомое изотонное продолжение $\mto{g}{X}{Y}$ отображения $f$ получаем как в~(\ref{s2}). Убедимся в существовании изотонного продолжения в случае $X$ имеющего произвольную  мощность. В соответствии с GCH мощность булеана $\mathfrak{B}(A)$ есть первый ординал строго больший $|A|$. Вложим $A$  в $\mathfrak{B}(A)$ с помощью $\mto{\nabla_A}{A}{\mathfrak{B}(A)}$ (см. замечание~\ref{r1}) и построим отображение $\mto{g_1}{\mathfrak{B}(A)}{Y}$, для которого $g_1\circ\nabla_A = f$, что можно сделать в соответствии с рассмотренным выше случаем. Рассматривая $\mathfrak{B}(A)$ как подмножество $\mathfrak{B}(X)$ и используя теорему~\ref{t5}, находим изотонное отображение $\mto{g_2}{\mathfrak{B}(X)}{Y}$, продолжающее $g_2$. В качестве искомого изотонного продолжения $\mto{g}{X}{Y}$ можно взять отображение $A\xrightarrow{\operatorname{in}_A}X\xrightarrow{\nabla_X}\mathfrak{B}(X)\xrightarrow{g_2}Y$, \[ \begin{diagram}
\node{\mathfrak{B}(A)}
\arrow[2]{e,t}{\operatorname{In}_A}
\arrow{se,t}{g_1}
\node[2]{\mathfrak{B}(X)}
\arrow{sw,t}{g_2}
\\
\node[2]{Y}
\\
\node{A}
\arrow[2]{n,l}{\nabla_A}
\arrow{ne,b}{f}
\arrow[2]{e,t}{\operatorname{in}_A}
\node[2]{X,}
\arrow[2]{n,r}{\nabla_X}
\end{diagram}\]
где $\operatorname{In}_A$ и $\operatorname{in}_A$~--- соответствующие вложения.

\textbf{(iii)$\Rightarrow$(i).} Пусть выполнено (iii), $A$ и $B$~--- подмножества $Y$,
\[ \max (\abs{A},\abs{B})<\alpha \]
и $a\ls{Y}b$ для любых $a\in A$ и $b\in B$. Достаточно доказать, что существует $y^*\in Y$ такое, что
\[ a\ls{Y}y^*\ls{Y}b \]
при всех $a\in A$ и $b\in B$. Предположим, что такого $y^*$ не существует. Тогда можно построить одноточечное расширение $\chu{\overline{Y}}$ ч.у. множества $\chu{Y}$ такое, что
\[ \overline{Y}=Y\cup{y^*}, \quad y^*\notin Y,\quad \ls{Y} \, = \, \ls{\overline{Y}}\cap(Y\times Y) \]
и
\begin{equation}\label{e4.9}
a\ls{\overline{Y}}y^*\ls{\overline{Y}}b
\end{equation}
при всех $a\in A$, $b\in B$. (Это легко сделать, используя вложение $Y\xrightarrow{\nabla_Y}\mfB (Y)$.) Положим
\[ X=A\cup B\cup\{y^*\}\, \text{ и  } \,  \ls{X} \, = \, \ls{\overline{Y}}\cap(X\times X). \]
По свойству (iii) вложение $\operatorname{in}\colon A\cup B\to Y$, $\operatorname{in}(t)=t$ для $t\in A\cup B$, имеет изотонное продолжение $\mto{g}{X}{Y}$. Из~(\ref{e4.9}) следует, что
\begin{equation}\label{e4.10}
a\ls{Y}g(y^*)\ls{Y}b
\end{equation}
при всех $a\in A$, $b\in B$. Так как $g(y^*)\in Y$, то~(\ref{e4.10}) противоречит сделанному выше предположению. Следовательно, $\chu{Y}$~--- $\alpha$-квазирешетка.

Импликация (iii)$\Rightarrow$(iv) очевидна, а проверка того, что имеет место (iv)$\Rightarrow$(ii) делается, как при доказательстве теоремы~\ref{t4.6*}.
\end{proof}

\begin{corollary}\label{c4.5*}
Пусть $\chu{Y}$~--- ч.у. множество. $\chu{Y}$ является полной решеткой тогда и только тогда, когда $\chu{Y}$~--- $\alpha$-квазирешетка для любого $\alpha$.
\end{corollary}
\begin{proof}
Следует из теоремы~\ref{t4.4*} и ~\ref{t4}.
\end{proof}

\begin{remark}\label{r*}
В стандартном определении $\beta$-универсальных множеств $\chu{P}$ (см, например, определение 2.1~\cite{Ha}) требуют вложимость в $P$ любого $\chu{X}$ с $\abs{X}=\beta$. Если $\alpha$~--- ординал, для которого $\beta=\aleph_{\alpha}$, то
\[ (|X|\leqslant \beta)\Leftrightarrow (|X|< \aleph_{\alpha+1}).\]
Следовательно, любое $\chu{P}$ являющееся $\beta$-универсальным  в стандартном смысле с $\beta =\aleph_{\alpha}$ будет $\aleph_{\alpha+1}$-универсальным в смысле определения~\ref{d4.3*}.

Из теоремы~\ref{t4.4*} легко вывести, что любое $\chu{P}$ являющееся $\beta$-универсальным в смысле определения~\ref{d4.3*} будет $\beta$-квазирешеткой. таким образом, универсальные ч.у. множества, исследовавшиеся ранее в известных работах~\cite{CD, Cue54, Cue58, Har68, Kur}  являются $\beta$-квазирешетками.
\end{remark}

\section{Полные локальные решетки и изотонные продолжения}
Пусть $\chu{X}$~--- ч.у. множество и $A$~--- подмножество $X$, имеющее наибольший и наименьший элементы, т.е. существуют $a^*$и $a_*$, принадлежащие $A$  такие, что
\begin{equation}\label{e5.1}
a_* \ls{X} a \ls{X} a^*
\end{equation}
для любого $a\in A$.  В этом случае
\begin{equation}\label{e5.2}
f(a_*) \ls{Y} f(a) \ls{Y} f(a^*)
\end{equation}
для любого изотонного $\mto{f}{A}{Y}$ и всех $a\in A$. Если $\mto{g}{X}{Y}$~--- изотонное продолжение $f$, то  будем говорить, что $g$ \emph{сохраняет экстремальные значения}, если
\begin{equation}\label{e5.3}
f(a_*) \ls{Y} g(x) \ls{Y} f(a^*)
\end{equation}
для всех $x \in X$.

\begin{definition}\label{d5.1}
Ч.у. множество $\chu{Y}$ будем называть \emph{полной локальной решеткой}, если интервал
\[ [y_*,y^*]_Y =\{ y\in Y, y_*\ls{Y}y\ls{Y}y^*\} \]
является полный решеткой для любых $y_*, y^*\in Y$, удовлетворяющих $y_*\ls{Y} y^*$.
\end{definition}
Очевидно, что любая полная решетка является полной локальной решеткой. Заметим также, что любая антицепь будет полной локальной решеткой. Таким образом, полная локальная решетка может не быть решеткой.

Лексикографическая сумма двух антицепей дает пример полной локальной решетки, не являющейся даже квазирешеткой, а квазирешетка, построенная в примере~\ref{ex7}, не является полной локальной решеткой. (Последнее следует из того, что не является  решеткой интервал $[-1,1]_R$, где  $R$~--- ч.у. множество, построенное в примере~\ref{ex7}.)

\begin{theorem}\label{t5.3}
Пусть $\chu{Y}$~--- непустое ч.у. множество. Следующие утверждения эквивалентны.
\begin{itemize}
  \item [(i)] $\chu{Y}$~--- полная локальная решетка.
  \item [(ii)] Для любого ч.у. множества $\chu{X}$ и любого имеющего наибольший и наименьший элементы множества $A\subseteq X$ любое изотонное отображение $\mto{f}{A}{Y}$ имеет изотонное продолжение на $X$, сохраняющее экстремальные значения.
\end{itemize}
\end{theorem}
\begin{proof}\textbf{(i)$\Rightarrow$(ii)} Пусть $\chu{Y}$~--- полная  локальная решетка, $\chu{X}$~--- ч.у. множество, $A\subseteq X$ и пусть $a^*$ и $a_*$~--- наибольший и, соответственно, наименьший элементы в $A$. Рассмотрим произвольное изотонное отображение $\mto{f}{A}{Y}$ и найдем изотонное продолжение $\mto{\Psi}{X}{Y}$ отображения $f$ такое, что ~(\ref{e5.3}) справедливо  для всех $x\in X$ с $g=\Psi$.

Так как $a^*$~--- наибольший элемент $A$, a $a_*$~--- наименьший элемент этого множества, то ~(\ref{e5.2}) имеет место для всех $a\in A$. Следовательно,
\[ f(A)\subseteq [f(a_*),f(a^*)]_Y. \]
В силу того, что $\chu{Y}$~--- полная локальная решетка, $[f(a_*),f(a^*)]_Y$~--- полная  решетка. Следовательно, по теореме~\ref{t4}, существует изотонное продолжение
\[ \mto{g}{X}{[f(a_*),f(a^*)]_Y} \]
отображения $A \ni a\mapsto f(a)\in [f(a_*), f(a^*)]_Y$. Теперь в качестве искомого $\mto{\Psi}{X}{Y}$ можно взять
\[ X\xrightarrow{g}[f(a_*),f(a^*)]_Y\xrightarrow{\operatorname{in}} Y, \]
где $\operatorname{in}(y)=y$ при всех $y\in [f(a_*),f(a^*)]_Y$.

\textbf{(ii)$\Rightarrow$(i).} Пусть выполнено (ii). Рассмотрим произвольный интервал $I=[y_*, y^*]_Y \subseteq Y$ и тождественное отображение $\mto{\operatorname{id}}{[y_*, y^*]_Y}{[y_*, y^*]_Y}$. Очевидно, что $y_*$  и $y^*$~--- экстремальные значения этого отображения. Пусть $\chu{L}$~--- полная решетка такая, что $L\supseteq I$ и $\ls{I}$ совпадает с сужением $\ls{L}$ на $I\times I$. Свойство (ii) влечет изотонную продолжаемость $\operatorname{id}$  на $L$. Применяя лемму~\ref{l3} находим, что $I$~--- полная решетка.
\end{proof}

Перейдем теперь к задаче построения изотонных продолжений для функций, заданных на $A$, имеющих наибольший и наименьший элементы и удовлетворяющих условию $\abs{A}<\alpha$.

Напомним, что подмножество $A$ ч.у. множества $X$ ограничено, если $\lc{A}\neq \varnothing \neq \uc{A}$.

\begin{definition}\label{d5.4}
Пусть $\alpha$~--- бесконечное кардинальное число. Ч.у. множество $\chu{Y}$ назовём \emph{локальной $\alpha$-квазирешеткой}, если для любых двух ограниченных $A,B \subseteq Y$, удовлетворяющих условиям~(\ref{e4.1}) и условию $\max(|A|,|B|)<\alpha$, найдется $y^*=y^*(A,B)$, для которого~(\ref{e4.2}) выполняется при всех $a\in A$ и $b \in B$.
\end{definition}

Очевидно, что любая $\alpha$-квазире\-шет\-ка является локальной $\alpha$-квази\-ре\-шет\-кой. Отметим также, что любая полная локальная решетка \mbox{$\chu{Y}$} есть локальная $\alpha$-квазирешетка при любом $\alpha$.  Действительно, пусть \linebreak[4] \mbox{$A, B \subseteq Y$} удовлетворяют условиям~(\ref{e4.1}) и условию $\max(|A|,|B|)<\alpha$ и пусть $b^*$~--- мажоранта $B$, а $a_*$~--- миноранта $A$. Тогда $A$ и $B$~--- подмножества полной решетки $I=[a_*,b^*]_Y$ и $y^*=\inf\nolimits_I B$ удовлетворяет~(\ref{e4.2}) для всех $a\in A$  и $b\in B$.

\begin{lemma}\label{l5.5}
Пусть $\alpha$~--- бесконечный кардинал, $\chu{Y}$~--- локальная $\alpha$-квази\-ре\-шет\-ка, $y_*, y^*\in Y$ и $y_*\ls{Y}y^*$. Тогда $I=[y^*,y_*]_Y$~--- $\alpha$-квази\-ре\-шет\-ка относительно порядка, индуцированного из $Y$.
\end{lemma}
\begin{proof}
Так как любые подмножества $A,B\subseteq I$  являются  ограниченными, то лемма следует непосредственно из определений~\ref{d4.2*} и~\ref{d5.4}.
\end{proof}

\begin{theorem}\label{t5.6}
Пусть $\alpha$~--- бесконечный кардинал, $\chu{Y}$~--- непустое ч.у. множество и $\chu{P}$~--- $\alpha$-универсально. Следующие утверждения эквивалентны.
\begin{itemize}
  \item [(i)] $\chu{Y}$~--- локальная $\alpha$-квазирешетка.
  \item [(ii)] Для любого ч.у. множества $\chu{X}$ и любого $A\subseteq X$, обладающего наибольшим и наименьшим элементами и удовлетворяющего условию $|A|<\alpha$, каждое изотонное отображение $\mto{f}{A}{Y}$ имеет изотонное продолжение, сохраняющее экстремальные значения.
  \item [(iii)] Для любого ч.у. множества $\chu{X}$ с $|X|<\alpha$ и любого, обладающего наибольшим и наименьшим элементами $A\subseteq X$, каждое изотонное отображение $\mto{f}{A}{Y}$ имеет изотонное продолжение, сохраняющее экстремальные значения.
  \item [(iv)] Для любого $A\subseteq P$, обладающего наибольшим и наименьшим элементами и такого, что $|A|<\alpha$, каждое изотонное отображение $\mto{f}{A}{Y}$ изотонно продолжается на $P$.
\end{itemize}
\end{theorem}

\begin{proof}
\textbf{(i)$\Rightarrow$(ii)}.  Пусть выполнено (i), $\chu{X}$~--- ч.у. множество  и пусть $A\subseteq X$ такое, что $|A|<\alpha$ и существуют $a_*, a^*\in A$, для которых
\[ a_* \ls{X} a \ls{X} a^* \]
при всех $a\in A$. Рассмотрим произвольное изотонное отображение $\mto{f}{A}{Y}$. Пусть $I:=[f(a_*),f(a^*)]_Y$. В силу леммы~\ref{l5.5} ч.у. множество $\chu{I}$ является $\alpha$-квазирешеткой. Следовательно,  по теореме~\ref{t4.4*},  отображение
\[ A\ni a\mapsto f(a) \in I \]
имеет изотонное продолжение $\mto{g}{X}{I}$. Очевидно, $g$ сохраняет экстремальные значения,
\[ f(a_*) \ls{Y} g(x) \ls{Y} f(a^*) \]
для всех $x\in X$. Отображение $X \xrightarrow{g} I \xrightarrow{\operatorname{in}} Y$ дает изотонное продолжение $f$, сохраняющее экстремальные значения.

Импликация (ii)$\Rightarrow$(iii) очевидна. Импликация (iii)$\Rightarrow$(i) может быть получена простой модификацией  рассуждений, приведенных при доказательстве соответствующей импликации в теореме~\ref{t4.4*}.
\end{proof}

\bigskip

{\bf А.\,А. Довгошей}

{\small  Институт прикладной математики и механики НАН Украины,}

{\small  ул. Розы Люксембург 74, Донецк, Украина,  83114,}

{\small e-mail: {\it aleksdov@mail.ru}
\end{document}